\documentclass[a4paper]{article}
\usepackage[french,english]{babel}  
\usepackage[T1]{fontenc}
\usepackage{amsmath,amssymb}  
\usepackage{color}
\usepackage{geometry}
\usepackage[all]{xy}
\usepackage{amsthm}
\usepackage{makeidx} 
\usepackage[nottoc]{tocbibind} 
\usepackage{graphicx}
\usepackage{varioref}
\usepackage{hyperref}
\usepackage{yhmath}

\hypersetup{
						colorlinks=true,
						breaklinks=true,
						urlcolor= blue,
						linkcolor= black,
						anchorcolor = black,
						citecolor = black,
						filecolor = black,
					}

\renewcommand\widering[1]{\ring{\wideparen{#1}}}

\newcommand{\wt}[1]{\widetilde{~#1}} 
			
\newcommand*{\longhookrightarrow}{\ensuremath{\lhook\joinrel\relbar\joinrel\longrightarrow}}

\newcommand{\C}{ {\mathbb C} }
\newcommand{\N}{ {\mathbb N} }
\newcommand{\Q}{ {\mathbb Q} }
\newcommand{\R}{ {\mathbb R} }
\newcommand{\Z}{ {\mathbb Z} }

\newcommand{\cB}{{\mathcal B}}
\newcommand{\cC}{{\mathcal C}}
\newcommand{\cE}{{\mathcal E}}

\newcommand{\cP}{{\mathcal P}}

\newcommand{\ds}{\displaystyle}

\newtheoremstyle{thm}
  {0pt}
  {5pt}
  {}
  {}
  {\bfseries}
  {:}
  {.5em}
  {}
  
\theoremstyle{thm}
\newtheorem{thm}{Theorem}
\labelformat{thm}{~#1}

\labelformat{cor}{~#1}
\newtheorem{prop}{Proposition}
\labelformat{prop}{~#1}

\newtheorem*{defn}{Definition}
\newtheorem*{demo}{Proof}
\newtheorem{exempl}{Example}
\labelformat{exempl}{~#1}

\newtheoremstyle{remarque}
  {3pt}
  {5pt}
  {}
  {}
  {\itshape}
  {:}
  {.5em}
  {}

\theoremstyle{remarque}
\newtheorem*{rmq}{Remark}

\date{}

\begin{document}

\title{Linear Gale transform of a starshaped spheres}
\author{Jerome Tambour}

\maketitle

\begin{abstract}
\noindent Gale transform is a simple but powerful tool in convex geometry. In particular, the use of Gale transform is the main argument in the classification of polytopes with few vertices. Many books and documents cover the definition of Gale transform and its main properties related to convex polytopes. But it seems that there does not exist document studying the Gale transform of more general objects, such that triangulation of spheres. In this paper, we study the properties of the Gale transform of a large class of such spheres called starshaped spheres.     
\end{abstract}

\tableofcontents

\section*{Introduction}

Gale transform is a very useful tool in the study of polytopes and other combinatorial objects. Mainly, it transforms a finite family $X$ of vectors in a vector space in another family $\overline{X}$ of vectors in a vector space of (generally) smaller dimension. Combinatorics of the first family is intimately related to the combinatorics of the second. For instance, if the vectors of $X$ are the vertices of a polytopes $P$, the combinatorial type of $P$ only depends on the combinatorics of $\overline{X}$ (cf. $\cite{G}$, paragraph $\S 5.4$). This allows a classification of polytopes with few vertices (a $d$-polytope is said to have few vertices if its number of vertices is smaller or equal to $d+3$) or to obtain useful informations on the classification on centrally-symmetric polytopes.
 
\vskip 5mm

Gale transforms also appears in toric geometry, i.e. the study of topological objects (mainly, algebraic varieties and smooth manifolds) endowed with an action of a compact torus $(S^1)^n$ or an algebraic torus $(\C^*)^n$. Main examples of those kind of "toric objects" are the well-known algebraic varieties (cf. $\cite{F}$, $\cite{Od}$ or $\cite{CLS}$) which play an important role in algebraic geometry. Many theorems of algebraic geometry has been shown first by studying the toric case and then proving the general case. More recently, many topological generalizations of toric varieties arose. The main and more fundamental examples are quasitoric manifolds (cf. $\cite{DJ}$), topological toric manifolds (cf. $\cite{FHM}$), torus manifolds (cf. $\cite{HM}$) and moment-angle complexs (cf. $\cite{BP}$). The common property of these objects is their combinatorial nature. Indeed, all of them can be described by combinatorial object and many of their properties of one object can be easily read on the associated combinatorial object. For instance, it is well-known that every toric varieties can be described in terms of a fan and that the variety is smooth (resp. compact, projective) if and only if the corresponding fan is regular (resp. complete, strongly polytopal). Topological toric manifolds are parametrized the so-called topological fans and torus manifolds by multifans.

\vskip 5mm

It appears that every toric variety and every quasitoric manifolds, provided it is an orbifold, can be obtained as the quotient of some moment-angle complex by a action of an compact torus. The family of moment-angle complexs is very large and is in one-to-one correspondence with the family of simplicial complexs. Many moment-angle complexs can have a very complicated topology (namely, a complicated structure for their cohomology ring). One of the main question concerning these complexs is the characterization of the moment-angle complexs which admits a complex structure. One recent result related to this question was obtained by the author in $\cite{T}$ (and also independently in $\cite{PU}$) and states that every moment-angle complex parametrized by a rationally starshaped sphere, that is simplicial sphere which are the underlying complex of a rational complete fan, can be endowed with a complex structure. Amongst those kind of complex manifolds, we recover all the complex structures on compact torus, Hopf and Calabi-Eckmann manifolds, some connected sums of products of spheres with odd dimension and intersections of quadrics (as studied in $\cite{LdM}$, $\cite{LdMV}$, or $\cite{M}$). Actually, the collection of moment-angle complexs associated to rationally starshaped spheres forms a sub collection of the family of LVMB manifold, constructed in $\cite{Bo}$ and in the three last articles above.

\vskip 5mm

More precisely, to construct LVMB manifolds, Lopez de Medrano, Verjovsky, Meersseman and Bosio studied an holomorphic action of $\C^m$ on some open subset in $\C^n$. The necessary and sufficient conditions for this action to be free, proper and cocompact (and, as a consequence, for the orbit space to be a compact complex manifold) are given in $\cite{Bo}$. We call these conditions \emph{Bosio conditions}. In $\cite{T}$, we show that, when Bosio conditions are satisfied, we can construct a rationally starshaped sphere (called associated complex) which encodes the topology of the orbit space. And conversely, for every such a sphere, we can find parameters for the holomorphic action satisfying Bosio conditions and such that the associated complex is the given sphere.

\vskip 5mm

The proof in $\cite{T}$ is geometric and uses the deep relation between LVMB manifolds and toric varieties. In this paper, we use Gale transform to show in a direct way the equivalence between Bosio conditions and the starshapeness of a sphere. The paper is divided in five sections. In the first section, we recall basic definitions of simplicial complexs and their realizations. The section $2$ is devoted to introduced the main object of study, namely the class of starshaped spheres. We also investigate some variations of the definition and prove in particular that every starshaped sphere admits a rational realization. In section $3$, we introduce and study two objects dual to starshaped spheres. We call them fundamental set and studiable systems. Actually, fundamentals sets are dual to pure simplicial complexes and in this section, we characterize fundamental sets which are dual to pseudomanifolds. Finally, in the forth section, we prove the main theorem of the paper and give a simple necessary and sufficient criterion to see if a realization of a simplicial complex is starshaped. This criterion use linear Gale transform and we hope it to be useful and lead to some interesting results. The fifth and last section contains some additional remarks and ideas for further investigations about linear Gale transforms.

\section{Notations and simplicial complexs}

\subsection{Some notations}

Firstly, we fix some notations which will be followed in this paper:

\begin{enumerate}
\item If $X=(x_j)_{j\in J}$ is a set whose index set is $J$ and $I$ is a subset of $J$, we note $X_I$ or $X(I)$ the set $(x_i)_{i\in I}$. 
\item If $J$ is a finite set, $|J|$ and $Card(J)$ are set for its cardinal.
\item If $A$ is a subset of a vector space $V$, then $Conv(A)$ (resp. $pos(A)$) is the convex hull (resp. the positive hull) of $A$. That means that 

\[
\ds{ Conv(A)=\left\{\ \sum_{j=1}^p \lambda_ja_j\ \middle/\ p\in\N,\ \lambda_j\in\R^+,\ a_j\in A,\ \sum_{j=1}^p\lambda_j=1\  \right\} }
\]

and

\[
\ds{ pos(A)=\left\{\ \sum_{j=1}^p \lambda_ja_j\ \middle/\ p\in\N,\ \lambda_j\in\R^+,\ a_j\in A\  \right\} }
\]

The relative interior of these sets will be $\widering{Conv}(A)$ et $\widering{pos}(A)$

\item Moreover, if $\cB=\left(e_1,\dots,e_p\right)$ is a basis of $E$, we note $\cB^*=\left(e_1^*,\dots,e_p^*\right)$ its dual basis in $E^*$. If $E=\R^n$ and $\cB$ is its canonical basis, then we identify $E$ and $E^*$ using the linear isomorphism $\phi$ defined by $\phi(e_j)=e_j^*$. 
\item If $E$ and $F$ are vector spaces and $L:E\longrightarrow F$ a linear map, then we note $L^*$ the linear map from $F^*$ into $E^*$ defined by $L^*(\phi)=\phi\circ L$. If $M$ is the matrix of $L$ in bases $\cB$ and $\cC$, then the matrix of $L^*$ in dual bases $\cC^*$ and $\cB^*$ is the transpose matrix of $M$.

\end{enumerate}

\subsection{Simplicial complexs}

We begin by recalling the definition of the combinatoric objects we will work with. Most of the definitions of this subsection are standards. We mainly follows $\cite{BP}$ for definitions and notations.

\vskip 5mm

\defn{Let $V$ be a finite nonempty set. An \emph{abstract simplicial complex}\footnote{or simply a \emph{simplicial complex}} on $V$ is a family $K$ of subsets of $V$ such that the two following properties are fulfilled:

\vskip 5mm

\begin{enumerate}
\item $\emptyset \in K$
\item $\forall \sigma\in K,\hskip 3mm \tau\subset\sigma \Rightarrow \tau\in K$. 
\end{enumerate}

\vskip 5mm

The elements of $K$ are \emph{faces} of $K$. Singletons in $K$ are called \emph{vertices} of $K$ and maximal elements of $K$ for inclusion are called \emph{facets}. Finally, the dimension of a face $\sigma$ of $K$ is $dim(\sigma)=|\sigma|-1$ and $dim(K)$ will be the maximum of the dimensions of its faces.
}

\vskip 5mm

\begin{defn} A \emph{geometric simplicial complex}\footnote{or simply a \emph{geometrical complex} } $C$ is a set of simplices\footnote{A simplex is the convex hull of a finite set of affinely independent points} in $\R^n$ which satisfies the following properties:

\vskip 5mm

\begin{enumerate}
\item If $\sigma$ is an element of $C$, then every facet (including the empty set) of $\sigma$ is an element of $C$.
\item If $\sigma$ and $\tau$ are elements of $C$, then $\sigma\cap\tau$ is a common face of $\sigma$ and $\tau$.
\end{enumerate}

\vskip 5mm

The \emph{support} of $C$ is the union $|C|$ of all its elements. We often identify $C$ to its support. 
\end{defn}

\vskip 5mm

Let $K$ be a simplicial complex of dimension $d$. A geometric complex $C$ is a \emph{realization} of $K$ if there exists a bijection between the vertex set of $K$ and the vertex set of $C$ such the image of a face of $K$ is the vertex set of a simplex of $C$. We can note that two realizations of the same simplicial complex are homeomorphic. \label{defsphere} We say that $K$ is a $d$-\emph{sphere} if $K$ admits a realization $C$ in $\R^q$ (with maybe $q$ not equal to $d+1$) whose support is homeomorphic to the $d$-sphere $S^d$. For instance, every polytopal complex is a simplicial sphere. Moreover, there exists simplicial spheres which are not polytopal (cf. $\cite{GS}$ or $\cite{Ba}$). We know that every simplicial sphere of dimension $2$ is polytopal (cf. $\cite{G}$ for example). Moreover, according to a theorem by Mani (cf. $\cite{Ma}$), every $d$-sphere with at most $d+4$ vertices is polytopal. In $\cite{GS}$, it is shown that there exists $39$ simplicial spheres with dimension $3$ and $8$ vertices (up to combinatorial equivalence), and $37$ of them are polytopal. The $2$ others are called Barnette sphere and Brückner sphere).

\section{Starshaped spheres} 
\label{defstars}

In this section, we introduce the main object of study of this paper. The starshaped spheres form an important family of triangulations of spheres because they are underlying fans of simplicial compact toric varieties. They also constitute a nice family of embeddable spheres (i.e. spheres of dimension $d$ which admits a realization in $\R^{d+1}$). In this section, beside stating the definition of these starshaped spheres, we prove that every such a sphere can be realized in a rational and starshaped way. We also discuss some variations around the definition. Firstly, let state the definition of starshapeness. We mainly follow $\cite{E}$.

\vskip 5mm  

\defn{Let $K$ be a $d$-dimensional simplicial complex. We say that $K$ is \emph{starshaped} if there exists a realization $|K|$ of $K$ in the euclidean space $\R^{d+1}$ and a point $p$ of $\R^{d+1}$ such that  every ray from $p$ intersects $|K|$ in only one point. We say that $|K|$ is \emph{starshaped} in $p$ and that $p$ is in the \emph{center} (or \emph{kernel}) of $|K|$. The center of $|K|$ is denoted by $Ker(|K|)$.
}

\vskip 5mm

\begin{figure}[h]
\centering
	\includegraphics[width=0.30\textwidth]{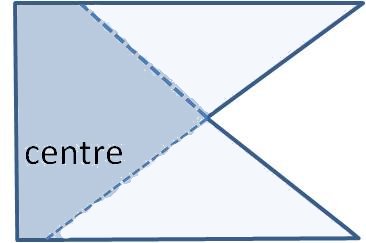}
\caption{Nonconvexe starshaped realization of the pentagon. The center is the darkest area in the interior of the pentagon.}	
\end{figure}

\vskip 5mm

\begin{rmq} The fact that a $d$-sphere admits a realization in $\R^{d+1}$ is an open problem (cf. \cite{MW}, \S 5).
\end{rmq}

\vskip 5mm

\begin{exempl} Every polytopal complex is a starshaped sphere. Indeed, every convex realization of a polytopal complex is starshaped and the center of the realization if exactly the (relative) interior of the realization. In $\cite{E2}$, Ewald and Schultz constructed a example of a nonstarshaped sphere (with dimension $3$ and $12$ vertices). 
\end{exempl}

\vskip 5mm

A starshaped complex is easily seen to be a simplicial sphere (see proposition \ref{starsphere}), so we will use indifferently the terms \emph{starshaped complex} or \emph{starshaped sphere}. We postpone the proof of this fact, in order to discuss some weakening of the definition.

\subsection{Weakening of the starshapedness}

There is two conditions in the definition of starshapedness of a simplicial complex $K$. The first one is that $K$ has to be embeddable, i.e. that it admits a realization in $\R^{dim(K)+1}$. The second one concerns the existence of a point such that every ray emanating from this point intersects the sphere in exactly one point. One may want to weaken one of these two conditions and studied the more general family of complexs obtained that way. Note that weakening the first condition implies to weaken the second one. Firstly, we prove the following proposition which shows that dropping the embeddability in the definition of a starshaped complex does not lead to an interesting notion.

\vskip 5mm

\begin{prop} Let $K$ be a simplicial complex on $\{1,\dots,n\}$ (with possibly ghost vertices). Let $(e_1,\dots,e_n)$ be the canonical basis of $\R^n$. The realization $|K|$ of $K$ in $\R^n$ whose simplices are exactly those of the form $Conv(e_j,j\in J)$ for $J\in K$, has the property that every ray from the origin intersect $|K|$ at at most one point.
\end{prop}

\vskip 5mm

\begin{demo} Let assume that there exists two nonzero colinear vectors $u$ and $v$ in the realization $|K|$. We denote $\sigma$ and $\tau$ the facets of $K$ containing $u$ and $v$ respectively. Then we can write $u=\sum_{j\in\sigma}t_je_j$ and $v=\sum_{j\in\tau}s_je_j$ for some positive numbers. But, since the vectors $\left(e_j\right)_{j=1}^n$ form a basis, considering $v-\lambda u=0$, we obtain that $\sigma=\tau$. This implies that the line $(u,v)$ is contained in $\text{Aff}(e_j,j\in\sigma)$. In particular, $0$ is an element of this affine hull. This is contradictory with the fact that the $e_j$'s are elements of a basis. $\square$
\end{demo}

\vskip 5mm

As a consequence, we make the following definition:

\vskip 5mm

\begin{defn} A simplicial complex $\cP$ of dimension $d$ which admits a realization $|\cP|$ in $\R^{d+1}$ such that there exists a point $p$ in $\R^{d+1}$ such that every ray emanating from $p$ intersects $|\cP|$ in at most one point is said \emph{weakly starshaped} (at $p$).
\end{defn} 

\vskip 5mm

\begin{prop}\label{starsphere} Let $\cP$ be a weakly starshaped complex. Then $\cP$ is the triangulation of a subset of a sphere. This is the whole sphere if and only if $\cP$ is starshaped.
\end{prop}

\vskip 5mm

\begin{demo} Since $\cP$ is weakly starshaped, then $\cP$ admits a realization $|\cP|$ in $\R^{d+1}$, such that any ray from the origin intersect $|\cP|$ in at most one point. Let $S$ be a sphere centered at the origin which does not intersect $\cP$ (this is possible since $|\cP|$ is compact). For every point $x$ of $|\cP|$, there is one point in $S$ on $\R_+x$. This defines an application $f$ from $|\cP|$ into $S$. This application is injective since $\cP$ is weakly starshaped. This prover the former part of the proposition. To prove the latter, it suffices to remark that the inverse of $f$ is the map which associate to $x$ in $f(|\cP|)\subset S$ the unique element of $\cP$ which is colinear to $x$. In particular, $\cP$ is a starshaped sphere if and only if $f^{-1}$ is defined on the whole sphere $S$, i.e. if and only if $\cP$ is a starshaped sphere. $\square$
\end{demo}

\subsection{Starshaped spheres and rationality}

Finding a center for a simplicial sphere is equivalent to solve some linear inequalities system. Indeed, if $K$ is a $d$-sphere realized in $\R^{d+1}$, then, according to Jordan theorem, its realization $|K|$ divides $\R^{d+1}\setminus |K|$ in two connected components: a bounded one, called "interior" of the sphere, and a unbounded component called "exterior" of the sphere. We can then make a coherent choice of normal vectors for the facets, choosing for example outwards normal vectors. 

\vskip 5mm

If $F$ is a facets of $|K|$, we note $F^o_+$ the open halfspace whose boundary is $\text{Aff}(F)$ and containing the outwards normal vectors of $F$ (the other halfspace being $F^o_-$). It is then clear that $|K|$ is a starshaped realization of $K$ if and only if $\displaystyle{\bigcap_{F facets}F^o_-}$ is non empty. In this case, $\displaystyle{\bigcap_{F facets}F^o_-}$ is exactly the center of the realization $|K|$.

\vskip 5mm

We now want to discuss the existence of a starshaped realization for a $d$-sphere whose vertices are all contained in some lattice\footnote{A lattice is the discrete group spanned by the element of a basis of $\R^{d+1}$} of $\R^{d+1}$. So, we make the following definition:

\vskip 5mm

\begin{defn} Let $K$ be a $d$-sphere. We say that $K$ is \emph{rationally starshaped} if there exists a lattice $L$ in $\R^{d+1}$, a point $p$ in $L$ and a realization $|K|$ of $K$ in $\R^{d+1}$ such that every vertex of $|K|$ is an element of $L$ and every ray emanating from $p$ intersects $|K|$ in exactly one point. We say that $|K|$ is \emph{rationally starshaped} with respect to $L$ and that $p$ is in the \emph{rational center} (with respect to $L$) of $|K|$. It will be denoted $Ker_L(|K|)$ in the sequel. Clearly, we have $Ker_L(|K|)=Ker(|K|)\cap L$.
\end{defn}

\vskip 5mm

In $\cite{T2}$, we prove that if $K$ is a rationally starshaped sphere, we can always assume that $L$ is the integer lattice $\Z^{d+1}$ and that $0$ is a center of $|K|$:

\vskip 5mm

\begin{prop} \label{reduction} Let $K$ be a rationally starshaped $d$-sphere. Then there exists a rationally starshaped realization $|K|$ of $K$ for the lattice $\Z^{d+1}$, with center $0$ and such that the canonical basis of $\Z^{d+1}$ belongs to the vertex set of $|K|$.
\end{prop}

\vskip 5mm

We can now prove the following theorem:

\vskip 5mm

\begin{thm} Every starshaped sphere is a rationally starshaped sphere.
\end{thm}

\vskip 5mm

\begin{demo} Let $K$ be a starshaped sphere. We note $d$ the dimension of $K$ and $v$ its number of vertices. 

\vskip 5mm

Let $p_1,\dots,p_v$ the different vertices of a starshaped realization $|K|$ of $K$ in $\R^{d+1}$. Up to a translation, we can assume that $0$ is in the center of $K$. By definition, if $F$ is a facet of $|K|$, then there exists a facet $J=\left(j_1,\dots,j_{d+1}\right)$ in $K$ such that $F=Conv\left(p_j, j\in J\right)$. We note $\Phi_J:\R^{d+1}\rightarrow \R$ the map defined by $\phi_J(x)=det\left(\wt{p_{j_1}},\dots,\wt{p_{j_{d+1}}},\wt{x}\right)$ and we have  

\[
x\in H_J=\text{Aff}(F) \Leftrightarrow \phi_J(x)=0
\]

We choose now an outwards normal vector $n_J$ for the facet $F$ and (up to change $\phi_J$ into $-\phi_J$) we can assume that $F^o_-=\{x\in \R^{d+1}\ / \phi_J(x)<0 \}$. In particular, $x$ is in the center of $|K|$ if and only if $\phi_J(x)<0$ for $J$ every $K$. 

\vskip 5mm

Actually, functions $\phi_J(x)$ depend upon $p=(p_1,\dots,p_v)$. So, we will note them $\phi_J(x,p)$ in the sequel. By construction, those functions are continuous. As a consequence, the map $\Phi_J:p\mapsto \phi_J(0,p)$ is continuous. But $\left(\R_{<0}\right)^{|J|}$ is open, so $U=\Phi_J^{-1}(\left(\R_{<0}\right)^{|J|})$ is an open set containing $p$. Since $\Q$ is dense in $\R$, we can find $\overline{p}=\left(\overline{p}_1,\dots,\overline{p}_v\right)$ in $U$ such that all $\overline{p}_j$ have rational coordinates. As a conclusion, up to a multiplication by an integer, then vectors of $\overline{p}$ are vertices of a rationally starshaped realization of $K$ (with center $0$). So $K$ is rationally starshaped. $\square$  
\end{demo}

\section{Fundamental sets and studiable systems}

In this section, we introduce the other important objects of this paper: fundamentals sets and studiable systems. In some sense, they are dual objects to starshaped spheres and their realizations. A fundamental set is the set of complements of facets of some pure simplicial complex. This generalizes the duality between simple and simplicial polytopes. A studiable system is obtained adding some geometric information, namely a family of vectors, to a fundamental set. They are particular of the torus graphs studied in $\cite{MMP}$.

\subsection{Fundamental sets}

In this subsection, $M$ is a fixed integer and $V$ a finite set whose cardinal $n$ is greater or equal to $M$. If the converse is not explicitly assumed, we will take $V=\{1,\dots,n\}$. 

\vskip 5mm

\begin{defn} A \emph{fundamental set} is a non empty set $\cE$ of subsets of $V$ such that every elements of $\cE$ has cardinal $M$. We will say that $\cE$ has type $(M,n)$. 
\end{defn}

\vskip 5mm

\begin{prop} The set of fundamental sets on $V$ with type $(M,n)$ is in bijection with the set of pure simplicial complexs of dimension $(n-M-1)$ on $V$. 
\end{prop}

\vskip 5mm

\begin{demo} A simplicial complex is exactly defined by its facets. If $K$ is a set of subsets of $V$, we note $Cl(K)$ for the unique simplicial complex whose facets are exactly the elements of $K$. The map which sends a fundamental set $\cE$ to the simplicial complex $Cl(\{\ V\setminus E\ /\ E\in\cE\ \})$ is easily seen to be bijective onto the set of pure simplicial complexs of dimension $(n-M-1)$. $\square$   
\end{demo}

\vskip 5mm

\begin{defn} Let $\cE$ be a fundamental set. The image of $\cE$ by the map in the previous proposition is called the \emph{associated complex} of $\cE$ and is denoted by $\cP(\cE)$ (or just $\cP$).  
\end{defn}

\vskip 5mm

Moreover, an element of $V$ is said to be \emph{indispensable} (with respect to $\cE$) if it is contained in every element of $\cE$. We will say that $\cE$ has type $(M,n,k)$ if it has type $(M,n)$ and has exactly $k$ indispensable elements. If $\cP$ is the simplicial complex associated to $\cE$, then $v\in V$ is an indispensable element of $\cE$ if and only if $v$ is a ghost vertex of $\cP$ (i.e. $\{v\}\notin\cP$). 

\vskip 5mm 

In the sequel, we will mainly discuss two combinatorial properties called respectively $SE$ (Substitute Existence) and $SEU$ (Substitute Existence and Uniqueness) properties:

\vskip 5mm

\[
(SE) \hskip 12mm \forall P\in {\cE},\ \forall k\in V,\hskip 2mm \exists\ k'\in P;\hskip 2mm (P\setminus\{k'\})\cup\{k\}\in{\cE}
\]

\vskip 5mm

\[
(SEU) \hskip 12mm \forall P\in {\cE},\ \forall k\in V,\hskip 2mm \exists!\ k'\in P;\hskip 2mm (P\setminus\{k'\})\cup\{k\}\in{\cE}
\] 

\vskip 5mm

We will also say that $\cE$ is \emph{minimal} for the $(SEU)$ property if it satisfies this property and if it does not strictly contain any fundamental set which satisfies the $(SEU)$ property.

\subsection{$(SEU)$ and pseudomanifold}
\label{sectpseudo}

In this subsection, we prove the following theorem:

\vskip 5mm 

\begin{thm}\label{minSEUpseudo} Let $\cE$ be a fundamental set of type $(M,n)$ with $n>M$, and $\cP$ its associated complex. Then $\cE$ is minimal for the $(SEU)$ property if and only if $\cP$ is a pseudomanifold. 
\end{thm}

\vskip 5mm

This theorem was already proved in $\cite{T}$. We include the proof here for sake of completeness. Firstly, we recall the definition of a pseudomanifold: 

\vskip 5mm

\begin{defn} \label{pseudo} Let $K$ be a simplicial complex. $K$ is a \emph{pseudo-manifold} if the two following properties are fulfilled:
\begin{enumerate}
	\item every codimension one face of $K$ is contained in exactly two facets.
	\item for all facets $\sigma,\tau$ of $K$, there exists a chains of facets $\sigma=\sigma_1,$ $\dots,$ $\sigma_r=\tau$ of $K$ such that $\sigma_i\cap\sigma_{i+1}$ is a codimension one face of $K$ for every $i\in\{1,\dots,r-1\}$.
\end{enumerate}
\end{defn}

\vskip 5mm

For instance, every simplicial sphere is a pseudo-manifold. More generally, a triangulation of a manifold (that is, a simplicial complex whose realization is homeomorphic to a topological manifold) is also a pseudo-manifold. The theorem above can be proved in two steps. The first step is to prove the equivalence between $(SEU)$ and the first property in the definition \ref{pseudo} above. The second consists in showing that a fundamental set satisfying $(SEU)$ is a disjoint union of fundamental set minimal for $(SEU)$. We now prove the first step.

\vskip 5mm

\begin{prop}[Step $1$] \label{crete} Let $\cE$ be a fundamental set. Then its associated complex $\cP$ satisfies the $SEU$ property if and only if every codimension one face of $\cP$ is contained in exactly two facets of $\cP$.
\end{prop}

\vskip 5mm

\begin{demo} Firstly, it is easy to restate the $(SEU)$ property for a fundamental set $\cE$ in terms of the associate complex $\cP$ of $\cE$. Actually, $\cE$ satisfies $(SEU)$ if and only if 

\[
\forall\ Q\in \cP_{max},\ \forall\ k\in \left\{1,\dots,n\right\},\ \exists!\ k'\notin Q;\ \left(Q\cup\{k'\}\right)\backslash\left\{k\right\}\in\cP_{max} 
\]

where $\cP_{max}$ is the set of facets of $\cP$.

\vskip 5mm

Now, we assume that $\cE$ verifies the $SEU$ property. Let $Q$ be a codimension one face of $\cP$. By definition, $Q$ is included in a facet $P$ of $\cP$. We put $P=Q\cup\{k\},\  k\in\{1,\dots,n\}\backslash Q$. By the remark above, there exists $k'\notin P$ (and so $k\neq k'$) such that $P'=(P\cup\{k'\})\backslash\{k\}$ is a facet of $\cP$. We have $P'=Q\cup\{k'\}$ so $Q$ is contained in at least two facets of $\cP$. Let assume that $Q$ is contained in a third facet $P''=Q\cup\{k''\}$. In this case, we have $P''=(P\cup\{k''\})\backslash\{k\}$, which contradicts the remark above.

\vskip 5mm

Conversely, let $Q$ be a facet of $\cP$ and $k\in\{1,\cdots,n\}$. If $k\in Q$, then $P=Q\backslash\{k\}$ is a codimension one face of $\cP$ and by hypothesis, $P$ is  contained in exactly two facets $Q_1$ and $Q_2$. One of them, say $Q_1$, is $Q$. The other is $Q_2$ and we have $Q_2=P\cup\{k'\}$. Then we have $k'\notin Q$ (on the contrary, we would have $Q_2=Q=Q_1$) and $Q_2=(Q\cup\{k'\})\backslash\{k\}$. Moreover, if $Q_3=(Q\cup\{k''\})\backslash\{k\}$ is another facet of $\cP$ with $k''\notin Q$, then $Q_3$ contains $P$ and by hypothesis, $Q_3=Q_2$ (i.e $k''=k'$). If $k\notin Q$, we remark that the element $k'\in\{1,\dots,n\}$ such that $Q'=(Q\cup\{k'\})\backslash\{k\}$ is a facet of $\cP$ is $k'=k$. Indeed, if $k'=k$, then $Q'=Q$ is a facet of $\cP$. And if $k'\neq k$, then $k\notin Q\cup\{k'\}$ and, as a consequence, $Q'=Q\cup\{k'\}$ is not in $\cP$.  $\square$
\end{demo}

\vskip 5mm

We can now prove the second step of the proof of the theorem:

\vskip 5mm

\begin{prop}[Step 2] Let $\cE$ be a fundamental set of type $(M,n)$ which verifies the $SEU$ property. Then, there exist an integer $p\in\N^*$, and fundamental sets $\cE_j$ of type $(M,n)$ which are minimal for the $SEU$ property and such that $\cE$ is the disjoint union $\displaystyle{\bigsqcup_{j=1}^p\cE_j}$.
\end{prop}

\vskip 5mm

\begin{demo} We proceed by induction on the cardinal of $\cE$. If $\cE$ is minimal for the $SEU$ property, then there is nothing to do. Let assume that it is not the case: there exists a proper subset $\cE_1$ of $\cE$ which is minimal for the $SEU$ property. We put $\overline\cE$ for its complement $\cE\backslash\cE_1$. It is obvious that $\overline\cE$ is a fundamental set (of type $(M,n)$). We claim that $\overline\cE$ verifies the $SEU$ property. Let $P$ be an element of $\overline\cE$ and $k\in\{1,\dots,n\}$. If $k\in P$, then, putting $k'=k$, we have that $(P\backslash\{k'\})\cup\{k\}=P$ is an element of $\overline\cE$. It is the only choice (for $k'$) since $P$ is in $\cE$ and $\cE$ verifies the $SEU$ property. Let assume now that $k$ is not in $P$. Since $P$ is an element of $\cE$, there exists exactly one $k'\in P$ such that $P'=(P\backslash \{k'\})\cup\{k\}$ is an element of $\cE$, too. We claim that $P'$ cannot be in $\cE_1$. Indeed, if it were the case, since $\cE_1$ is minimal for the $SEU$ property, there would exist exactly one $k''\in P'$ such that $P''=(P'\backslash\{k''\})\cup\{k'\}\in\cE_1$. But $P=(P'\backslash\{k\})\cup\{k'\}$ is in $\cE$ and $k\in P'$. So, $k''=k$ and $P''=P$. As a consequence, $\overline\cE$ is a fundamental set of type $(M,n)$ which verifies the $SEU$ property with cardinal strictly smaller than $\cE$. Applying the induction hypothesis on $\overline\cE$, we have the  decomposition of $\cE$ we were looking for. $\square$
\end{demo}

\vskip 5mm

Finally, we can use these two steps to prove the theorem \ref{minSEUpseudo}. To prove it, we introduce the following graph:

\vskip 5mm

\begin{defn} Let $\cE$ be a fundamental set of type $(M,n)$. We define the (unoriented) graph $\Gamma$ by requiring that its vertices are fundamental subsets of $\cE$ and two vertices $P$ and $Q$ are related by an edge  if and only if there exist $k\notin P,k'\in P$ such that $Q=(P\backslash\{k'\})\cup\{k\}$. Equivalently, we relate two subsets of $\cE$ if and only if they differ exactly by one element. $\Gamma$ is called the \emph{replacement graph} of $\cE$.
\end{defn}

\vskip 5mm

\begin{exempl} If $\cP$ corresponds to the boundary of a simplicial polytope $P$, then $\Gamma$ is the nerve of the dual of $P^*$. 
\end{exempl}

\vskip 5mm

\begin{rmq} It is clear that if $\cE$ satisfies the $(SEU)$ property, then $\Gamma$ is a regular $(n-M-1)$-valent graph. Moreover, each minimal fundamental set in the partition of $\cE$ in the step $2$ correspond to a connected component of $\Gamma$. In particular, $\cE$ is minimal for the $(SEU)$ property if and only if $\Gamma$ is connected. 
\end{rmq}

\vskip 5mm

\begin{demo}[proof of the theorem \ref{minSEUpseudo}] Let assume that $\cE$ is minimal for the $SEU$ property. This implies that every codimension one face of $\cP$ is contained in exactly two facets. Now, let $\sigma,\tau$ be two distinct facets of $\cP$. So, $P=\sigma^c$ and $Q=\tau^c$ are two elements of $\cE$. By minimality for the $SEU$ property, $\Gamma$ is connected. Consequently, there exists a sequence $P_0=P,P_1,\dots ,P_r=Q$ of fundamental subsets such that $P_i$ and $P_{i-1}$ differ by exactly one element. We denote $R_i$ the set $P_{i-1}\cup P_i$ with $M+1$ elements. Its complement $R_i^c$ is thus a codimension one face of $\cP$. If we put $\sigma_i=P_i^c$, we have $R_i^c=\sigma_{i-1}\cap\sigma_i$ so $\sigma_0=\sigma,\dots ,\sigma_r=\tau$ is the sequence connecting the facets $\sigma$ an $\tau$. Consequently, $\cP$ is a pseudo-manifold. The converse is proved in a similar way. $\square$
\end{demo}

\subsection{Studiable systems}

\defn{A \emph{studiable system} of type $(M,n)$ is a couple $(\cE,\Lambda)$ where $\cE$ is a fundamental set with type $(M,n)$ and $\Lambda=\left(\lambda_1,\dots,\lambda_n\right)$ a family of vectors in $\R^M$ such that for every $P$ in $\cE$, the vectors of $(\lambda_p)_{p \in P}$ span $\R^M$ as a vector space.}

\vskip 5mm

\defn{Let $(\cE,\lambda)$ be a studiable system. We say that this system fulfills the \emph{imbrication condition} $(Imb)$ if, for every elements $P$ and $Q$ in $\cE$, the (relative) interiors of $pos(\lambda_p,p\in P)$ and of $pos(\lambda_q,q\in Q)$ are non disjoint.}

\vskip 5mm

\begin{rmq} The imbrication condition $(Imb)$ is equivalent to the following:

\vskip 5mm

\[ 
\begin{array}{rcl}
\widetilde{(Imb)} & & \forall P,Q\in\cE,\ \widering{Conv}(\{0\}\sqcup\{\lambda_p,p\in P\})\cap\widering{Conv}(\{0\}\sqcup\{\lambda_q,q\in Q\})\neq\emptyset   
\end{array}
\]

\end{rmq}

\vskip 5mm

\begin{prop}\label{ouvertetudiable} Let $\cE$ be a fundamental set with type $(M,n)$. Then the set of families $\lambda=\left(\lambda_1,\dots,\lambda_n\right)$ of $n$ vectors in $\R^{M}$ such that $(\cE,\lambda)$ is studiable is an open set $S$ of $\left(\R^M\right)^n$.
\end{prop}

\vskip 5mm

\begin{demo} For every $P=\{p_1<\dots<p_M\}$ in $\cE$, we define $f_P:\left(\R^M\right)^n\rightarrow\R$ by 

\[
f_P(\lambda)=det(\wt{\lambda_{p_1}},\dots,\wt{\lambda_{p_M}})
\]

\vskip 5mm

and $f:\left(\R^M\right)^n\rightarrow\R^{|\cE|}$ by $f(\lambda)=\left(f_P\left(\lambda\right)\right)_{P\in \cE}$. Then $S=f^{-1}((\R^*)^{|\cE|})$. But $f$ is continuous (since it is a polynomial function of the coefficients of the $\lambda_j$) and $(\R^*)^{|\cE|}$ is open, so $S$ is open. $\square$
\end{demo}

\section{Linear Gale transform of a starshaped sphere}

\subsection{Linear Gale transform}

\label{defgale}

We recall the definition of a linear Gale transform and some of its basic properties. Our main source is the section $\S 4$ of the Chapter II of \cite{E}. Proof without statement in this subsection can be found in this book or are left to the reader. Basically, linear Gale transform generalizes the construction of the dual basis of a basis in a finite dimensional vector space. This transform induces a duality between a family of vectors in a vector space and another family of vectors in a vector space with (generally) smaller dimension. Properties a family can be translated into properties of the other family. 

\vskip 5mm

Let $E$ be a vector space and let consider $X=(x_1,\dots,x_n)$, a family of elements of $E$. We assume that $X$ linearly spans $E$. We fix now a vector space $F$ with basis $(b_1,\dots,b_n)$ and we define a surjective linear map $L_1$ from $F$ onto $E$ setting $f(b_j)=x_j$ for every $j$. Finally, let $G$ be a vector space and $L_1$ a linear map from $G$ into $F$ such that the sequence

\[
0 \longrightarrow G \stackrel{L_2}{\longrightarrow} F \stackrel{L_1}{\longrightarrow} E \longrightarrow 0
\]

is exact \footnote{$G$ and $L_2$ always exist. We can choose $W=Ker(f)$ and $L_2$ to be the inclusion $G \longhookrightarrow F$ for instance.}. 

\vskip 5mm
 
Then, basic results of linear algebra imply that the dual sequence 
 
\[
0 \longleftarrow G^* \stackrel{g^*}{\longleftarrow} F^* \stackrel{f^*}{\longleftarrow} E^* \longleftarrow 0
\] 

is also exact.

\vskip 5mm

\begin{defn} The family $\overline{X}$ of vectors $(\overline{x}_1,\dots,\overline{x}_n)$ defined by $\overline{x}_j=g^*(b_j^*)$ for every $j$ is a \emph{linear Gale transform} of $X$. A linear Gale transform is defined up to the choice of $G$ and $L_2$, so defined up to linear isomorphism.
\end{defn} 

\vskip 5mm

\begin{exempl} Let $X=\left(x_1,\dots,x_r,e_1,e_2,\dots,e_{d+1}\right)$ a family of vectors in $E=\R^{d+1}$. We note $A$ the $(d+1,r+d+1)$-matrix whose columns are the vectors $x_j$ and $x^1,\dots,x^{d+1}$ the rows of $A$. Then the family $\overline{X}=(e_1,\dots,e_r,-x^1,\dots,-x^{d+1})$ is a linear Gale transform of $X$.
\end{exempl}

\vskip 5mm 

The most important property of linear Gale transform is that it transforms linearly independent vectors in a $d$-dimensional space into spanning vectors of a $(n-d)$-dimensional space. More precisely, we have:

\vskip 5mm

\begin{prop}[$\cite{E}$]\label{indep} Let $X=(x_1,\dots,x_n)$ a family of vectors which span a vector space $E$ and $\overline{X}$ a linear Gale transform of $X$ in some space $G^*$. Let $I$ be a subset of $V=\{1,\dots,n\}$. Then:
\begin{enumerate}
 \item  The family $X(I)$ is linearly independent in $E$ if and only if $\overline{X}(V\setminus I)$ spans $G^*$.
 \item  In particular, $X(I)$ is a basis of $E$ if and only if $\overline{X}(V\setminus I)$ is a basis of $G^*$.
\end{enumerate}
\end{prop}

\vskip 5mm

Linear Gale transforms and its affine version (cf. $\cite{E}$) have proved their usefulness in studying the combinatorics of polytopes with few vertices. Indeed, if $X=\{x_1,\dots,x_n\}$ is a family of vectors spanning a $d$-dimensional vector space, then its linear Gale transform $\overline{X}$ lies in a $(n-d)$-dimensional vector space. So, when $n$ is not too big compared to $d$, it is often simpler to study the combinatorics of $\overline{X}$ than the combinatorics of $X$. One of the main steps in this study is a criterion based on Gale transform to know if a family $X$ is the set of vertices of a polytope. In the end of next section, we give a general criterion to recognize if a family of vectors are vertices of a realization of a starshaped sphere. We call this criterion \emph{Bosio conditions}.

\subsection{A bunch of properties}

Let $\cP$ be a pure simplicial complex on $\{1,\dots,n\}$ with dimension $d$ which admits a realization $|\cP|$ in $\R^{d+1}$. Let $x_1,\dots,x_n$ be the vertices of $|\cP|$ and $\cP_{max}$ the set of facets of $\cP$. We put $\cE=\{\ P^c\ /\ P\in \cP_{max} \ \}$. We also note $\overline{X}=\left(\overline{x}_1,\dots,\overline{x}_n\right)$ a linear Gale transform of $X=\left(x_1,\dots,x_n\right)$.

\vskip 5mm

\textbf{Notation:} If $P$ is an element of $\cP$, we put $C_P=Conv\left(x_p,p\in P\right)$ and $\sigma_P=pos\left(x_p,p\in P\right)$. We also note $H_P=\text{Aff}\left(x_p,p\in P\right)$. Conversely, if $E$ is an element of $\cE$, we put $D_E=Conv\left(\overline{x}_j,j\in E\right)$, and  $\delta_E=pos\left(\overline{x}_j,j\in E\right)$.

\vskip 5mm

Clearly, by definition of a realization of a simplicial complex, we have 

\[
|\cP|=\left\{\ C_P\ \middle/\ P\in\cP\ \right\}
\]

\vskip 5mm

Moreover, if $P\subset \{1,\dots,n\}$, we set $P_0=\{0\}\sqcup P$. Finally, we note $\cE_0=\{\ E_0\ /\ E\in\cE\ \}$ and we put $\overline{x_0}=0$.

\vskip 5mm

\begin{defn} We define some properties for the couple $(\cE,\overline{X})$: 

\[
\begin{array}{rcl}
 (gen)             & & \forall E\in\cE,\ \left(\overline{x_j},j\in E\right)\ \text{is a basis of}\ \R^M 										        \\
 \widetilde{(Imb)} & & \forall P,Q\in\cE,\ \widering{D}_{P_0}\cap\widering{D}_{Q_0}\neq\emptyset          										      \\
 (Imb)             & & \forall P,Q\in\cE,\ \widering{\delta}_{P}\cap\widering{\delta}_{Q}\neq\emptyset   								      			\\
\end{array}
\]

\end{defn}

\vskip 5mm

\begin{defn} One the other side, at the level of the $\cP$ and $X$, we define natural properties for a simplicial complex of dimension $d$:

\[
\begin{array}{rcl}
 (simpl) & & \forall P\in\cP_{max},\ \left(x_p,p\in P\right)\ \text{is a basis of}\ \R^{d+1}  \\
 (Sep)   & & \forall P\neq Q\in\cP,\ \widering{\sigma}_{P}\cap\widering{\sigma}_{Q}=\emptyset \\
 (wStar) & & |\cP| \text{ is weakly starshaped in } 0                                                \\  
\end{array} 
\]

\end{defn}

\vskip 5mm

In the sequel, our main task will be to study the relationship between these latter properties are equivalent to the former. Firstly, the proposition $4.11$ of $\cite{E}$ means that $(gen)$ is equivalent to $(simpl)$. It is also obvious that $(Imb)$ is equivalent to $\widetilde{(Imb)}$.

\subsection{Study of the starshapedness}

To begin, we can show the following usual property of starshapedness:

\vskip 5mm

\begin{prop} $(wStar)$ implies $(simpl)$ and $(Sep)$.
\end{prop}

\vskip 5mm

\begin{demo} Let assume there exists a facet $P$ of $\cP$ such that $\left(x_p,p\in\cP\right)$ is not a basis of $\R^{d+1}$. Since $\cP$ has dimension $d$, this means that $\left(x_p,p\in P\right)$ does not span $\R^{d+1}$. As a consequence, there exists a linear hyperplane $H$ containing $\left(x_p,p\in P\right)$. Then $C_P$ is included in $H$. But $0$ is also an element of $H$, so every semi-line from $0$ passing through a point $x\neq 0$ in the interior of $C_P$ intersects $C_P$ (and $|\cP|$) in more than one point. Then $|\cP|$ is not weakly starshaped. 

\vskip 5mm

Let $P$ and $Q$ be two facets of $\cP$. We assume that there exists an element $x$ in $\widering{\sigma}_P\cap\widering{\sigma}_Q$. Then there exists  $\lambda_P>0$ and $\lambda_Q>0$ such that $\lambda_Px\in \widering{C}_P$ and $\lambda_Qx\in \widering{C}_Q$. But $x$ is non zero and $|\cP|$ is weakly starshaped in $0$, so $\lambda_Px=\lambda_Qx$. We deduce from this that $\lambda_P=\lambda_Q$. But $C_P\cap C_Q$ is a face of $|\cP|$. It follows that there exists $R\in\cP$ such that $C_P\cap C_Q=C_R$. We have $x\in C_R$ and $C_R$ is a face of $C_P$. If $C_R$ is a proper face of $C_P$, we have $C_R\cap\widering{C}_P=\emptyset$. This is a contradictory, so $C_R=C_P$. With the same reasoning, we obtain $C_R=C_Q$. As a consequence, we have $P=Q$.
$\square$
\end{demo}

\vskip 5mm

And the converse holds:

\vskip 5mm 

\begin{prop} $(Sep)$ and $(simpl)$ implies $(wStar)$.
\end{prop}

\vskip 5mm

\begin{demo} Let assume firstly that there exists $P\in\cP$ such that $0\in H_P$. Then $\left(x_p,p\in P\right)$ does not span $\R^{d+1}$, and as a consequence $(simpl)$ is not satisfied. Let assume now that $|\cP|$ is not weakly starshaped in $0$ but that $(simpl)$ is satisfied. Let show that in this case, $(Sep)$ is not fulfilled. Since $|\cP|$ is not weakly starshaped in $0$, then there exists a point $x\in\R^{d+1}\setminus\{0\}$ such that $\R_{\geq0}x\cap |\cP|$ contains at least two elements $u$ and $v$. The property $(simpl)$ is fulfilled, so there exists two distinct elements $P$ and $Q$ in $\cP$ such that $u\in\widering{C}_P$ and $v\in\widering{C}_Q$ (if the contrary holds, then there would exist a linear hyperplane such $C_P$, the facet containing $u$ and $v$). But $u$ and $v$ are nonzero elements of $\R_{\geq0}x$, so there exists $\lambda>0$ such that $v=\lambda u$. As a consequence, $v$ is an element of $\widering{\sigma}_P\cap\widering{\sigma}_Q$. Choosing facets $\overline{P}$ and $\overline{Q}$ in $\cP$ containing $P$ and $Q$ respectively, we obtain $\widering{\sigma}_{\overline{P}}\cap\widering{\sigma}_{\overline{Q}}\neq\emptyset$. So $(Sep)$ is not satisfied. $\square$  
\end{demo}

\vskip 5mm

This proves that, for a simplicial complex, the weakly starshapedness is equivalent of being \emph{fan-like} (cf. $\cite{E}$ definition $5.6$ and remark below p.$91$). We now prove an important theorem, which give a criterion for a realization of a simplicial complex to be weakly starshaped, based on properties on its linear Gale transform. 

\vskip 5mm

\begin{thm}\label{th1} $(gen)$ and $(Imb)$ are equivalent to $(simpl)$ and $(Sep)$.
\end{thm}

\vskip 5mm

\begin{demo} Let assume that $\cP$ fulfills $(gen)$ and $(Imb)$ but that $|\cP|$ is not weakly starshaped in $0$. Then there exists two nonzero vectors $u$ and $v$ which are positively colinear (i.e. there exists $\lambda>0$ such that $v=\lambda u$) and two distinct facets $P$ and $Q$ in $\cP$ such that $u\in\widering{C}_P$ and $v\in\widering{C}_Q$. Then $P^c$ and $Q^c$ are elements of $\cE$. Since $(Imb)$ holds, we have $\widering{\delta}_{P^c}\cap\widering{\delta}_{Q^c}\neq\emptyset$. Let consider $\alpha$ an element of $\widering{\delta}_{P^c}\cap\widering{\delta}_{Q^c}$. Then there exists $(t_p)_{p\notin P}$ and $(s_q)_{q\notin Q}$ two families of strictly positives numbers such that 

\vskip 5mm

\[
\ds{\alpha=\sum_{p\notin P}t_p\overline{x}_p=\sum_{q\notin Q}s_q\overline{x}_q }
\] 

\vskip 5mm

We then put $t_p=0$ for every $p\in P$ and similarly $s_q=0$ for every $q\in Q$. And then we put $r_j=t_j-s_j$ for every $j$. Hence we obtain that 

\vskip 5mm

\[
\ds{\sum_{j=1}^nr_j\overline{x}_j=0}
\]

\vskip 5mm

But $\overline{x}_j=L_2(b_j^*)$ (using notation of section \ref{defgale}). So we have 

\vskip 5mm

\[
\ds{L_2^*(\sum_{j=1}^nr_jb_j^*)=0}
\]

\vskip 5mm

i.e. $\ds{\psi=\sum_{j=1}^nr_jb_j^*}$ is in $Ker(L_2^*)=Im(L_1^*)$. As a consequence, there exists a linear form $\phi$ on $\R^{d+1}$ such that $L_1^*(\phi)=\phi\circ L_1=\psi$. We denote the hyperplane $Ker(\phi)$ by $H$. Since we have $L_1(b_j)=x_j$ for every $j$, we obtain that $\phi(x_j)=r_j$ for every $j$.

\vskip 5mm

But we know that

\vskip 5mm

\[
r_j=\left\{
\begin{array}{rl}
 t_j>0,\ & \text{if}\ j\in Q\setminus P \\
 -s_j<0,\ & \text{if}\ j\in P\setminus Q \\
 0,\ & \text{if}\ j\in P\cap Q 		  \\
 ?,\ & \text{if}\ j\in P^c\cap Q^c 
\end{array}
\right.
\]

\vskip 5mm

So $C_P$ is included in the half space $H_+=\{\ x\in\R^{d+1} /\ \phi(x)\geq 0\ \}$ and $u$ belongs to its interior $H^o_+$. Similarly, $C_Q$ is included to the other half space $H_-$ and $v$ to its interior $H_-^o$. But $0$ is an element of $H$, so the open ray emanating from the origin $0$ and passing through $u$ (and $v$), that is the set $\R_{>0}u$ is contained in only one of the open halfspace $H^o_+$ or $H^o_-$. This is absurd, so $|\cP|$ is starshaped in $0$.

\vskip 5mm

Conversely, let assume that the conditions $(simpl)$ and $(Sep)$ are satisfied. This implies that $(gen)$ is satisfied. Let show that $(Imb)$ is also satisfied. For, we assume this is not the case and show that it would lead to a contradiction. So, there exists two set $E$ and $F$ in $\cE$ such that $\widering{\delta}_{E}\cap\widering{\delta}_{F}$ is empty. We set $P=E^c$ and $Q=F^c$ and $P$ and $Q$ are facets of $\cP$. Because of $(Sep)$, we have $\widering{\delta}_P\cap\widering{\delta}_Q=\emptyset$. So there exists a linear form $\phi$ on $\R^{d+1}$ such that, if $H=Ker(\phi)$, we have $x_j\in H$ whenever $j\in P\cap Q$, $x_p\in H^o_+$ if $p\in P\setminus Q$ and $x_q\in H^o_-$ if $q\in Q\setminus P$ if $q\in Q\setminus P$. 

\vskip 5mm

We then put $r_j=\phi(x_j)$ and $\ds{\psi=\sum_{j=1}^nr_j \overline{x}_j=0}$. Remark also that $r_j=0$ whenever $j\in P\cap Q$. So we have

\vskip 5mm

\[
\ds{0= \sum_{j\in P\cap Q}r_j\overline{x}_j+\sum_{j\in P\setminus Q}r_j\overline{x}_j+\sum_{j\in Q\setminus P}r_j\overline{x}_j+\sum_{j\in P^c\cap Q^c}r_j\overline{x}_j}
\]

\vskip 5mm  

Then, we put $r_j^+=Max(r_j,0)\geq 0$ and $r_j^-=Max(0,-r_j)\geq 0$. We also choose an arbitrary strictly positive number $\lambda$. Then we have
 
\[
\ds{
\sum_{j\in P^c\cap Q^c}\left(r_j^-+\lambda\right)\overline{x}_j+\sum_{j\in P\setminus Q}\left(-r_j\right)\overline{x}_j=\sum_{j\in Q\setminus P}r_j\overline{x}_j+\sum_{j\in P^c\cap Q^c}\left(r_j^++\lambda\right)\overline{x}_j
}
\]

\vskip 5mm

Finally, since $Q^c=\left(P^c\cap Q^c\right)\sqcup\left(P\setminus Q\right)$, we obtain that $\ds{ \alpha=
\sum_{j\in P^c\cap Q^c}\left(r_j^-+\lambda\right)\overline{x_j}  +  \sum_{j\in P\setminus Q} \left(-r_j\right)\overline{x_j} }$ is an element of $\widering{\delta}_{P^c}\cap\widering{\delta}_{Q^c}$. This contradicts a previous assumption. So $(Imb)$ is satisfied. $\square$
\end{demo}

\subsection{Relation between $(SE)$ and imbrication condition}

In the previous subsection, we proved that weakly starshapedness of a simplicial complex is equivalent to the imbrication condition $(Imb)$ for its linear Gale transform. But we also proved previously that if $\cP$ is starshaped, then its realization $|\cP|$ is a sphere, hence a pseudomanifold. This implies that $\cE$ and $\cE_0$ are minimal for the $(SEU)$ property. So, we proved the following proposition:

\vskip 5mm

\begin{prop} Let $\cP$ be a starshaped sphere and $X=(x_1,\dots,x_n)$ the set of vertices of a starshaped realization at the origin of $\cP$. Then we denote $\cE$ the set of complements of facets of $\cP$ and $\overline{X}$ a linear Gale transform of $X$. Then $\left(\cE,\overline{X}\right)$ satisfies $(gen)$, $(Imb)$ and the minimality for the $(SEU)$ property.
\end{prop}
 
\vskip 5mm

The main theorem of this paper is to prove the converse, i.e. that if $(\cE,\Lambda)$ is a studiable system satisfying $(Imb)$ and minimal for $(SEU)$, then its associate complex is starshaped. As we show in the sequel of this subsection, only the Substitute Existence $(SE)$ property is needed beside imbrication condition $(Imb)$ to guarantee the starshapedness of the associate complex. We call the conditions $(Imb)$ and $(SE)$ for a studiable system the \emph{Bosio conditions}. Indeed, we will state that:

\vskip 5mm 

\begin{thm} Bosio conditions $(SE)$, $(Imb)$ and $(gen)$ imply $(SEU)_{min}$. 
\end{thm}

\vskip 5mm

Actually, this fact was already know in $\cite{B}$. The easy part of the proof is to show that the three conditions $(SE)$, $(Imb)$, and $(gen)$ imply the property $(SEU)$. It can be done only using affine geometry. The following proof come from $\cite{B}$; we rewrite it here for the sake of completeness:

\vskip 5mm

\begin{prop} A fundamental system satisfying $(SE)$, $(Imb)$ and $(gen)$ also satisfies $(SEU)$.
\end{prop}

\vskip 5mm

\begin{demo} Let $(\cE,\Lambda)$ be a fundamental system which satisfies the three conditions above but that $\cE$ does not satisfies the uniqueness of the substitute. This means that there exists an element $P$ in $\cE$ and $k$ in $\{1,\dots,n\}$ such that there exists at least two elements $k'$ and $k''$ in $P$ such that $Q'=\left(P\setminus\{k'\}\right)\cup\{k\}$ and $Q''=\left(P\setminus\{k''\}\right)\cup\{k\}$ are also in $\cE$. 

\vskip 5mm

We consider then $Q=P\setminus\{k',k''\}$. If $H$ denote the linear subspace spanned by the vectors $\lambda_j$, with $j\in Q\sqcup\{k\}$, then $H$ is an linear hyperplane (because of $(gen)$) and  the imbrication condition $(Imb)$ forces $\lambda_{k'}$ and $\lambda_{k''}$ to be in the same open halfspace delimited by $H$. But it is clear (using coordinates for example) that the hyperplane $H'$ defined by $(\lambda_q)_{q\in Q\sqcup\{k'\}}$ or the hyperplane $H''$ defined by $(\lambda_q)_{q\in Q\sqcup\{k''\}}$ strictly separates $\lambda_k$ and $\lambda_{k''}$ or strictly separates $\lambda_k$ and $\lambda_{k'}$. This contradicts $(Imb)$. $\square$
\end{demo}

\vskip 5mm

The more difficult part of the proof is to show that Bosio conditions are sufficient to guarantee the minimality for the $(SEU)$ property. This is also proved in $\cite{B}$, and the beautiful and short proof there uses topological and complex geometric arguments (and the relation between Bosio conditions and complex manifolds called LVMB manifolds). We give here another and longer proof of this fact, using only affine geometry. 

\vskip 5mm

\begin{prop} A fundamental system satisfying $(SE)$, $(Imb)$ and $(gen)$ is minimal for the $(SEU)$ property.
\end{prop}

\vskip 5mm

\begin{demo} Let $(\cE,\Lambda)$ be a fundamental system satisfying $(SE)$, $(Imb)$ and $(gen)$ and whose type is $(M,n)$. Then we consider a family $X=\left(x_1,\dots,x_n\right)$ a vectors in $\R^{d+1}$ whose linear Gale transform is $\Lambda=\left(\lambda_1,\dots,\lambda_n\right)$  (actually, $X$ is itself a linear Gale transform of $\Lambda$). We also defined $\cP$ to be the associated complex of $\cE$, so its facets are the complements of the elements of $\cE$. As a consequence, $\cP$ is a disjoint union of pseudomanifolds and the vectors of $X$ constitute the vertices of a weakly starshaped realization of $\cP$. Let $\Sigma$ be the fan whose generators are the vectors of $X$ and whose underlying complex is $\cP$. If $\Sigma$ is a complete fan, then $\cP$ is starshaped and $\cP$ is a sphere. In particular, $\cP$ is a pseudomanifold, and this is equivalent for $\cE$ to be minimal for the $(SEU)$ property.

\vskip 5mm

Let assume that $\Sigma$ is not complete and let consider a nonzero vector which is not in the support of $\Sigma$. Since there is only a finite number of cones, the minimum of the distances between $v$ and the maximal cones of $\Sigma$ is obtained for some cone $\sigma$ of $\Sigma$ and at some point $w$ in $\sigma$. Because the property $(simpl)$ is satisfied, $\sigma$ has dimension $d+1$ (where $d$ is the dimension of $\cP$), so $w$ is not an interior point of $\sigma$. Slightly moving $v$, we can assume that $w$ belongs to some facets $\tau$ of $\sigma$. It is clear that the linear hyperplane $H$ defined by $\tau$ separates $\sigma$ and $v$ (otherwise, the distance between $v$ and $\sigma$ would be obtained in other facets of $\sigma$). Let call $H^+$ the halfspace delimited by $H$ and containing $\sigma$ and $H^-$ the one containing $v$. Let remark that $v$ is the interior of $H^-$.  

\vskip 5mm

Finally, the facets $\tau$ corresponds to a codimension $1$ face $Q$ of $\cP$. But $Q$ has to be contained into two facets of $\cP$, since $\cP$ is a disjoint union of pseudomanifolds. This is equivalent to say that $\tau$ has to be a facets of two maximal cones of $\Sigma$, which one is $\sigma$. The other maximal cone, say $\nu$, has to be contained in $H^-$. But this implies that the distance between $v$ and $\sigma$ is not the minimum of the distance between $v$ and the maximal cones of $\Sigma$. 

\vskip 5mm

As a consequence, we obtain that $\Sigma$ is complete and thus $\cE$ is minimal for the $(SEU)$ property. $\square$
\end{demo}

\vskip 5mm

To conclude this subsection, we can remark that there exists examples of fundamental sets $\cE$ satisfying minimality for the $(SEU)$ property but such that there exists no family $\Lambda$ of vectors of $\R^M$ such that $(\cE,\Lambda)$ is a studiable system satisfying the imbrication condition $(imb)$. Examples were already given in $\cite{B}$. We can here give a method to create many examples: take any pseudo-manifold $\cP$ which is not a triangulation of a part of a sphere and take $\cE$ to be its set of complements of facets. Then $\cE$ is minimal for the $(SEU)$. Moreover, if $(\cE,\Lambda)$ is a studiable system satisfying the imbrication condition $(imb)$, then $|\cP|$ is a triangulation of a subset of a sphere, which is a contradiction to our previous assumption.  

\vskip 5mm

Finally, we are in position to show the following theorem:

\vskip 5mm

\begin{thm}\label{thmprinc} Let $\cP$ be a simplicial complex and $|\cP|$ a realization of $\cP$ whose set of vertices is $X$. We denote $\cE$ the set of complements of the facets of $\cP$ and $\Lambda$ a linear Gale transform of $X$. Then, $|\cP|$ is starshaped at the origin if and only if $(\cE,\Lambda)$ is a studiable system satisfying Bosio conditions $(Imb)$ and $(SEU_{min})$. 
\end{thm}

\vskip 5mm

To prove the theorem, we only have to show the following proposition:

\vskip 5mm

\begin{prop} Let $(\cE,\Lambda)$ be a studiable system satisfying Bosio conditions and $\cP$ its associated complex. Then $|\cP|$ is starshaped at $0$.  
\end{prop}

\vskip 5mm

\begin{demo} First of all, the conditions $(gen)$ and $(imb)$ imply that $|\cP|$ is a triangulation of a subset $A$ of a sphere. In particular, $A$ is closed and has a non empty boundary if $A$ is not the whole sphere. But $(imb)$ and $(SE)$ imply that $\cE$ satisfies minimality for the $(SEU)$ property, so $\cP$ is a pseudomanifold. As a consequence, $|\cP|$ has no boundary and thus $A$ is the whole sphere. This justifies that $|\cP|$ is a starshaped realization of $\cP$. $\square$
\end{demo}

\section{Additional remarks and conclusions}

\subsection{Condition $(K)$}

In the study of $LVMB$ manifolds, where the notion of Bosio conditions was introduced, an additional condition is sometimes really useful to investigate studiable systems (and the manifolds described by these systems). We call this additional condition \emph{condition (K)}:

\vskip 5mm

\defn{Let $(\cE,\lambda)$ be a studiable system with type $(M,n)$. We say that this system satisfies the $(K)$ condition if there exists an automorphism of $\R^M$ such that, for every $j$, $\phi(\lambda_j)$ has coordinates in $\Z^M$.
}

\vskip 5mm

Let $(\cE,\Lambda)$ be a studiable system satisfying Bosio condition and $\cP$ its associated complex. By \ref{thmprinc}, $\cP$ is a starshaped sphere. In section \ref{defstars}, we show that every starshaped sphere admits a rationally starshaped realization,say $|\cP|$ (whose vertices are non necesarly the elements of a linear Gale transform of $\Lambda$). Let denote $X$ the set of vertices of $|\cP|$ and $\Lambda'$ a linear Gale transform of $X$. Then, the elements are $\Lambda'$ have coordinates in $\Q^M$. Then, multiplying by a well-chosen integer and using the continuity of the Gale transform and the density of $\Q$ in $\R$, we obtain:

\vskip 5mm

\begin{prop} Let $\cE$ be a fundamental set satisfying $(SE)$. Then the set of families $\Lambda$ of vectors of $\R^M$ such that $(\cE,\Lambda)$ is a studiable system satisfying $(Imb)$, $(SE)$ and the condition $(K)$ is dense in the set of families $\Lambda$ of vectors of $\R^M$ such that $(\cE,\Lambda)$ is a studiable system satisfying $(Imb)$ and $(SE)$. 
\end{prop}

\vskip 5mm

Actually, this can be proved directly without using linear Gale transform and the starshapedness of the associated complex. This follows from the following property:

\vskip 5mm 

\begin{prop} Let $\cE$ be a fundamental set with type $(M,n)$. Then the set $\tilde{S}$ of families $\lambda=\left(\lambda_1,\dots,\lambda_n\right)$ of $n$ vectors in $\R^M$ such that $(\cE,\lambda)$ is studiable and fulfills the imbrication condition $(Imb)$ is an open set of $S$ (where $S$ is the open set of the proposition \ref{ouvertetudiable}) and of $\left(\R^M\right)^n$.
\end{prop}

\vskip 5mm

We can also prove this fact directly: 

\vskip 5mm

\begin{demo} The proof is elementary but needs some notations. Firstly, let fix $P=\{p_1<\dots<p_{2m+1}\}$ an element of $\cE$. And let $j$ be an element of $\{1,\dots,2m+1\}$. Then:

\vskip 5mm

\begin{enumerate}
\item $H^P_j$ the affine hyperplane in $\C^m=\R^{2m}$ spanned by the vectors $\left(\lambda_{p_k}\right)_{k\neq j}$.
\item Then there exists an affine form $\Phi^P_j$ such that $H^P_j=Ker(\Phi^P_j)$. We choose $\Phi^P_j$ in such a way that $\lambda_{p_j}$ is a element of the open halfplane $(H^P_j)^o_+$. In particular, we can put

\[
\phi_j^P(x)=det\left(\wt{\lambda_{p_1}},\dots,\wt{\lambda_{p_{j-1}}},\wt{x},\wt{\lambda_{p_{j+1}}},\dots, \wt{\lambda_{p_{2m+1}}}\right)
\]

and

\[ 
\Phi_j^P=\frac{|\phi_j^P(\lambda_{p_j})|}{\phi_j^P(\lambda_{p_j})}\phi_j^P
\]

We note also $L(\Psi_j^P)$ the linear form associated to $\Phi_j^P$.

\item We note $n_j^P$ the outwards normal vector to $H_j^P$ defined by $n_j^P=\left(L(\Psi_j^P)(e_1),\dots,L(\Psi_j^P)(e_{2m})\right)$. This vector is constructed in such a way that it is directed into $(H^P_j)^o_+$ and we have 

\[
\Psi_j^P(x)=<n_j^P,x>
\]

(for the usual scalar product of $\R^{2m}$). 
\item We put $\alpha_j^P$ to be the isobarycenter of vectors $\left(\lambda_{p_k}\right)_{k\neq j}$. In particular, $\alpha_j^P$ is an element of $H^P_j$. And if $\alpha$ is in $\R^{2m}$, we put 

\[
\Psi_j^P(\alpha)=<\alpha-\alpha_j^P,n_j^P>
\]

We then have the relation 

\[
\alpha\in (H^P_j)^o_+ \Leftrightarrow \Psi_j^P(\alpha)>0
\]

\item We also note $\Psi^P(\alpha)=(\Psi_1^P(\alpha),\dots,\Psi_{2m+1}^P(\alpha))$ and obviously we have 

\[
\alpha\in \widering{Conv}(\lambda_p,p\in P) \Leftrightarrow \Psi^P(\alpha)\in (\R_{>0})^{2m+1}
\]
\end{enumerate}

\vskip 5mm

Now, if $P$ and $Q$ are two elements of $\cE$, we note $\Theta_{P,Q}(\alpha)=\left(\Psi^P(\alpha),\Psi^Q(\alpha)\right)$. The relative interiors of $Conv(\lambda_p,p\in P)$ and of $Conv(\lambda_q,q\in Q)$ have a common point if and only if $\Theta_{P,Q}(\alpha)$ is an element of $\left(\R_{>0}\right)^{2(2m+1)}$.

\vskip 5mm

Finally, if $\alpha=(\alpha_{P,Q})_{P,Q\in\cE}$ is an element of $(\R^{2m})^{|\cE|^2}$, we define $\Theta(\alpha)=\left(\Theta_{P,Q}(\alpha_{P,Q})\right)_{P,Q\in\cE}$. Then $(\cE,\lambda)$ satisfies the imbrication condition if and only if there exists $\alpha=(\alpha_{P,Q})_{P,Q\in\cE}$ in $(\R^{2m})^{|\cE|^2}$ such that $\Theta(\alpha)$ is an element of $\left(\R_{>0}\right)^{2(2m+1)|\cE|^2}$.  

\vskip 5mm

All the previous objects depend actually on $\lambda$. In the sequel, we will note $\Theta(\alpha;\lambda)$ for $\Theta(\alpha)$ (and in a similar way for other objects). The determinant is continuous, so the functions $\phi_j^P(x,\lambda)$ are all continuous. Similarly, the isobarycenter map $\alpha_j(\lambda)$ is continuous, so the map $\Theta(\alpha,\lambda)$ is also continuous.

\vskip 5mm

We can finally prove the proposition. Let assume that $(\cE,\lambda)$ is a studiable system fulfilling the imbrication condition. Then, what is above implies that $\Theta$ is well defined and that there exists $\alpha=\left(\alpha_{P,Q}\right)_{P,Q\in\cE}$ such that $\mu=\Theta(\alpha,\lambda)$ is an element of $\left(\R_{>0}\right)^{2(2m+1)|\cE|^2}$. This last set is an open set so there exists an open neighborhood $V$ of $\mu$ in $\left(\R_{>0}\right)^{2(2m+1)|\cE|^2}$. But $\Theta$ is continuous so $\lambda\mapsto \Theta(\alpha,\lambda)$ also. Particularly, $\Theta(\alpha,\cdot)^{-1}(V)$ is an open set of $S$. As $S$ is itself open in $\left(\R^{2m}\right)^n$, $\tilde{S}$ is actually an open set of $\left(\R^{2m}\right)^n$. $\square$
\end{demo}

\subsection{Conclusion}

As we tried to show, linear Gale transform can be used to characterize combinatorial objects more general than vertex set of polytopes. In this paper, we gave a criterion to show if a family of vectors are the vertices of a starshaped, or equivalently, generators of a complete simplicial fans. It could be interesting to try to characterize through linear Gale transform other objects which appears in toric geometry such as characterist

\bibliographystyle{short}
\bibliography{Bibliographie}

\end{document}